\documentclass[10pt]{article}

\usepackage{amsmath}
\usepackage{amsfonts}
\usepackage{amssymb}

\begin{document}

\title{ Rational Approximation, Hardy Space - Decomposition of Functions in $L_p, p<1$: Further Results in Relation to Fourier Spectrum Characterization of Hardy Spaces}
\author{Guantie Deng and Tao Qian}

\author{Guantie Deng \thanks{School of  Mathematical Sciences, Beijing
Normal University, Beijing, 100875. Email:denggt@bnu.edu.cn.\ This
work was partially supported by NSFC (Grant  11271045) and by SRFDP
(Grant 20100003110004)} \ \ and
        Tao Qian \thanks{Department of Mathematics, University of Macau, Macao (Via Hong
Kong). Email: fsttq@umac.mo. The work was partially supported by
Multi-Year Research Grant (MYRG) MYRG116(Y1-L3)-FST13-QT, Macao Science and Technology Fund FDCT 098/2012/A3}}

\date{}
\maketitle
\begin{abstract}
    Subsequent to our recent work on Fourier spectrum characterization of Hardy spaces $H^p(\mathbb{R})$ for the index range $1\leq p\leq \infty,$ in this paper we prove further results on rational Approximation, integral representation and Fourier spectrum characterization of functions in the Hardy
  spaces $H^p(\mathbb{R}), 0 < p\leq \infty,$ with particular interest in the index range $ 0< p \leq 1.$  We show that the set of  rational functions in $ H^p(\mathbb{C}_{+1}) $ with the single pole $-i$
  is dense in $ H^p(\mathbb{C}_{+1}) $ for $0<p<\infty.$  Secondly, for $0<p<1$,  through rational function approximation we show that any function $f$ in $L^p(\mathbb{R})$ can be decomposed  into a sum  $g+h$,
   where $g$ and $h$ are, in the $L^p(\mathbb{R})$ convergence sense, the non-tangential boundary
limits of functions in, respectively, $ H^p(\mathbb{C}_{+1})$ and
$H^{p}(\mathbb{C}_{-1}),$ where
$H^p(\mathbb{C}_k)\ (k=\pm 1) $  are the Hardy spaces in
the half plane  $ \mathbb{C}_k=\{z=x+iy: ky>0\}$. We give  Laplace integral representation formulas for functions in the Hardy spaces $H^p,$ $0<p\leq2.$ Besides one in the integral representation formula we give an alternative version of Fourier spectrum characterization for functions in the boundary Hardy spaces $H^p$ for $0<p\leq 1.$
\end{abstract}


\section{Introduction}

  The classical Hardy space $ H^p(\mathbb{C}_k) ,\ 0<p< +\infty, k=\pm 1 $,
consists of the functions $f$ analytic in the
half plane $ \mathbb{C}_k=\{z=x+iy: k y>0\}.$ They are Banach spaces for $1\leq p<\infty$ under the norms
$$
 \|f\|_{H^p_k}=\sup_{k y>0}\left
(\int_{-\infty}^{\infty}|f(x+iy)|^pdx\right )^{\frac{1}{p}};
    $$
    and complete metric spaces for $0<p<1$ under the metric functions
    $$ d(f,g)=\sup_{k y>0}\int_{-\infty}^{\infty}|f(x+iy)|^pdx.$$
 A function  $f\in  H^p(\mathbb{C}_k)$ has non-tangential
boundary limits ({\it NTBL}s) $f(x)$ for almost all $x\in \mathbb{R}.$ The corresponding
boundary function belongs to $L^p(\mathbb{R}).$ For $1\leq p<\infty,$
$$
\|f\|_p=\left
(\int_{-\infty}^{\infty}|f(x)|^pdx\right )^{\frac{1}{p}}=\|f\|_{H^p_k}.
$$
For $p=\infty$ the Hardy spaces $H^\infty (\mathbb{C}_k) \ ( k=\pm 1)$ are defined to be the set of bounded analytic functions in $\mathbb{C}_k.$ They are Banach spaces under the norms
 \[ \| f\|_{H^\infty _k}={\rm sup}\{ |f(z)|\ :\ z\in \mathbb{C}_k\}.\]
 As for the finite indices $p$ cases any $f\in H^\infty (\mathbb{C}_k)$ has non-tangential boundary limit (NTBL) $f(x)$ for almost all $x\in \mathbb{R}.$  Similarly, we have
 \[ \|f\|_\infty ={\rm ess\ sup}\{ |f(x)|\ :\ x\in \mathbb{R}\}=\|f\|_{H^\infty (\mathbb{C}_k)}.\]

 We note that  $g(z)\in  H^p(\mathbb{C}_{-1})$ if and only if
the function $f(z)=\overline{g(\bar{z})} \in  H^p(\mathbb{C}_{+1})$.
The correspondence between their non-tangential
boundary limits and the functions themselves in the Hardy spaces is an isometric isomorphism. We denote by
$ H^p_k(\mathbb{R}) $ the spaces of the non-tangential boundary limits, or, precisely,
$$H^{p}_{k}
(\mathbb{R})=\bigg\{f\ : \mathbb{R}\to \mathbb{C},  f \mbox { is the NTBL of a function in }\ H^p(\mathbb{C}_k)\bigg\}.
$$
For $p=2$ the Boundary Hardy spaces $H^2_k(\mathbb{R})$ are Hilbert spaces.

We will need some very smooth classes of analytic functions that are dense in $ H^p(\mathbb{C}_{+1})$ and will play the role  of the polynomials in the disc case. J,B. Garnett in \cite{Ga} shows the following results.\\

 {\bf  Theorem A} (\cite{Ga}) Let N be a positive integer. For $0<p<\infty, \ pN>1$, the class $\mathfrak{A}_N$ is dense in  $ H^p(\mathbb{C}_{+1})$, where  $ \mathfrak{A}_N$ is the family of $ H^p(\mathbb{C}_{+1})$ functions satisfying

 (i) $f(z)$ is infinitely differentiable in $ \overline{\mathbb{C}_{+1}}$,

 (ii) $|z|^Nf(z)\rightarrow 0 $ as $z\rightarrow \rightarrow \infty , z\in \overline{\mathbb{C}_{+1}}$.\\

We shall notice that the condition $pN>1$ implies  the class $\mathfrak{A}_N$ is contained  in  $ H^p(\mathbb{C}_{+1})$. Let $ \alpha
$ be a complex number and $\mathfrak{R}_N(\alpha )$  the family of rational functions $f(z)=(z+\alpha)^{-N-1}P((z+\alpha )^{-1})$, $P(w)$ are  polynomials. We  notice that  the class $\mathfrak{R}_N(\alpha)$ is contained  in  the class $\mathfrak{A}_N$ for
$ \mbox{Im}\alpha >0$.

The tasks of this paper are three-fold. The first, replacing  the class $\mathfrak{A}_N $  by the class $\mathfrak{R}_N(i)$,  we will generalize Theorem A  as

  {\bf  Theorem 1}  Let N be a positive integer. For $0<p<\infty, \ Np>1$, the class $\mathfrak{R}_N(i)$ is dense in  $ H^p(\mathbb{C}_{+1})$.

  {\bf  Corollary 1}  Let N be a positive integer. For $0<p<\infty, \ Np>1$, the class $\mathfrak{R}_N(-i)$ is dense in  $ H^p(\mathbb{C}_{-1})$.

The second task is decomposition of functions in $L^p(\mathbb{R}), 0<p<1,$ into sums of the corresponding Hardy space functions in $H^p_{+1}(\mathbb{R})$ and in $H^p_{-1}(\mathbb{R})$ through rational functions approximation, and, in fact, by using what we call as rational atoms.

{\bf  Theorem   2 }\ (Hardy Spaces Decomposition of $L^p$ Functions For $0<p<1$) \ Suppose that   $0<p<1$ and  $ f\in L^p(\mathbb{R})$.
Then,  there exist a positive constant $A_p$ and  two  sequences of rational functions $\{P_k(z)\}$ and $\{Q_k(z)\}$ such that  $ P_k\in H^p(\mathbb{C}_{+1})$, $ Q_k\in
H^p(\mathbb{C}_{-1})$
and
\begin{eqnarray}\label{cite2}
\sum_{k=1}^{\infty} \left (\|P_k\|^p_{H^p_{+1}}+\|Q_k\|^p_{H^p_{-1}}\right )\leq A_p\|f\|_p^p ,
\end{eqnarray}
\begin{eqnarray}\label{cite3}
   \lim_{n\rightarrow \infty}||f-\sum_{k=1}^{n}(P_k+Q_k)||_p
  =0.
\end{eqnarray}
Moreover,
\begin{eqnarray}\label{cite4}
g(z)=\sum_{k=1}^{\infty}P_k(z)\in H^p(\mathbb{C}_{+1}) ,\ \ h(z)=\sum_{k=1}^{\infty}Q_k(z)\in H^p(\mathbb{C}_{-1}),
\end{eqnarray}
and
 $g(x)$ and $ h(x)$ are the non-tangential boundary
values   of  functions for $g\in H^p(\mathbb{C}_{+1})$ and $h\in
H^p(\mathbb{C}_{-1})$, respectively, $f(x)=g(x)+h(x)$ almost
everywhere, and
$$\|f\|_p^p\leq  \|g\|_p^p+\|h\|_p^p\leq A_p\|f\|_p^p , $$
  that is, in the sense of $L^p(\mathbb{R})$,
$$
L^p(\mathbb{R})=H^{p}_{+1}
(\mathbb{R})+H^p_{-1}(\mathbb{R}).
$$

  For the uniqueness of the decomposition, we can ask the following question: what is  the intersection space $H^p_{+1}(\mathbb{R})\bigcap H^p_{-1}(\mathbb{R})$? A.B. Aleksandrov (\cite{Ale} and \cite{Cim} )  gives an answer for this problem.

{\bf Theorem B} (\cite{Ale} and \cite{Cim} )
Let $0<p<1$ and ${X}^p$ denote the $L^p$ closure of the set of $f\in L^p(\mathbb{R})$ which can be written in the form
$$f(x)=\sum\limits_{j=1}^N\frac{c_j}{x-a_j},\ \ a_j\in \mathbb{R},\ c_j\in \mathbb{C}.$$
Then $${X}^p= H_{+1}^p(\mathbb{R})\bigcap  H_{-1}^p(\mathbb{R}).$$

 A.B. Aleksandrov's proof  (\cite{Ale} and \cite{Cim} ) is rather long involving vanishing moments and the Hilbert transformation. We present a more straightforward proof for this result.

The Fourier transform of a function $f\in L^1({\Bbb R})$  is defined, for $ x\in{\Bbb R}$, by
$$
\hat{f}(x)=\frac{1}{\sqrt{2\pi}}\int_{{\Bbb R}}f(t)e^{- ix t}dt.
$$

Based on the Fourier transformation defined for $L^1(\mathbb{R})$-functions, Fourier transformation can be extended to $L^2(\mathbb{R}),$ and then to $L^p(\mathbb{R}), 0<p<2,$ and finally to $L^p (\mathbb{R}), 2<p\leq \infty,$ the latest bring in the distribution sense.

The classical Paley-Wiener Theorem deals with the Hardy $H^2(\mathbb{C}_{+1}) $ space (\cite{Den},\cite{Du}, \cite{Ga},\cite{Koo} and \cite{Ru}) asserting that
$f\in L^2(\mathbb{R})$ is the NTBL of a function in $H^2(\mathbb{C}_{+1})$ if and only if ${\rm supp}\hat{f}\subset [0,\infty).$
  Moreover, in such case,
   the integral representation
\begin{eqnarray}\label{PW}
f(z)=\frac{1}{\sqrt{2\pi}}\int_0^\infty e^{itz}\hat{f}(t)dt\end{eqnarray}
holds.

We recall that Fourier transform of a tempered distribution $T$ is defined through the relation
$$
(\widehat{T}, \varphi)=(T,\hat{\varphi}) \ \
$$
for $\varphi$ in the Schwarz class $\mathbb{S}.$ This coincides with the traditional definition of Fourier transformation for functions in  $L^p(\mathbb{R}),1\le p\leq 2$. A measurable function $f$ satisfying
$$
\frac{f(x)}{(1+x^2)^m}\in L^p(\mathbb{R}) \ \  ( 1\leq p\leq \infty)
$$
for some positive integer integer $m$ is called a tempered $L^p$ function (when $p=\infty$ such a function is often called a slowly increasing function). The Fourier transform is a one to one mapping from $\mathbb{S}$ onto $\mathbb{S}$.

It is proved in \cite{Q1} that a function in $H^p_{+1}(\mathbb{R}), 1\leq p\leq \infty,$ induces a tempered distribution $T_f$ such that ${\rm supp}\hat{T}_f\subset [0,\infty).$ In \cite{Q2}, the converse of the result is proved: Let $T_f$ be the tempered distribution induced by  $f$ in $L^p(\mathbb{R}), 1\leq p\leq \infty.$ If ${\rm supp}\hat{T}_f\subset [0,\infty),$ then $f\in H^p_{+1}(\mathbb{R}).$

 The third task of this paper is to extend the above mentioned Fourier spectrum results, as well as the formula (\ref{PW}) to $0< p<1$.

 {\bf  Theorem 3 }\ (Integral Representation Formula For Index Range $0<p\leq 1$) \ If  $0<p\leq 1,\ f\in H^p(\mathbb{C}_{+1})$, then
there exist a positive constant $A_p$, depending only on $p$,  and
 a slowly increasing  continuous function $F$ whose support is contained in $[0,\infty)$,  satisfying that, for $\varphi$ in the Schwarz class $\mathbb{S}$,
$$
(F, \varphi)=\lim_{y>0, y\rightarrow0}\int_{\mathbb{R}}f(x+iy)\hat{\varphi}(x)dx, \ \
$$
and that
\begin{eqnarray}\label{cite5}
|F(t)|\leqslant A_p\|f\|_{H_{+1}^p}|t|^{\frac{1}{p}-1}, \ \ \ (t\in \mathbb{R})
\end{eqnarray}
 and
 \begin{eqnarray}\label{cite6}
 f(z)=\frac{1}{\sqrt{2\pi}}\int^\infty_0 F(t)e^{itz} dt \quad\quad (z\in \mathbb{C}_{+1}).
\end{eqnarray}

 P. Duren cites on page 197 of \cite{Du} that the argument to prove the integral representation (\ref{PW}) for $p=2$
 can be generalized to give an analogous representation for $H^p(\mathbb{C}_{+1})$-functions for
$1\leq p<2$. A proof for the range $1\leq p<2,$ in fact, is not obvious, and so far has not appeared in the literature, as far as concerned by the authors. We are to prove the following theorem corresponding to what Duren stated.

{\bf  Theorem C }\ (\cite {Du}, integral Representation Formula For Index Range $1\leq p \leq 2$) \ Suppose  $1\leq p\leq 2,\ f\in L^p(\mathbb{R})$. Then $f\in  H^p_{+1}(\mathbb{R})$ if and only if  supp$\hat{f}\subset[0,+\infty)$.
Moreover, under such conditions the integral representation (\ref{PW}) holds.

 We, in fact, prove analogous formulas for all the cases $0<p\leq 2.$ For the range $0<p<1$ we need to prove extra estimates to guarantee the integrability (See the proof of Theorem 3). The idea of using rational approximation is motivated by the studies of Takenaka-Malmquist systems in Hardy $H^p$ spaces for $1\leq p\leq \infty$ (\cite{Wa}, \cite{QW}).  For the range of $1\leq p\leq \infty$ this aspect is related to the Plemelj formula in terms of Hilbert transform that has immediate implication to Fourier spectrum characterization in the case. For the range of $0<p<1$ the Plemelj formula approach is not available.

\section{Proofs of theorems }

We need the following Lemmas.

{\bf  Lemma  1 }\ Suppose that   $0<p<1$ and  $R$ is a rational
function with $ R\in L^p(\mathbb{R}).$  For $k=\pm 1$, if  $R(z)$ is analytic in the half plane $ \mathbb{C}_{k}$, then
$ R\in H^p(\mathbb{C}_k)$.

{\bf Proof }\ \ \  Let  $0<p<1,\ R(z)=\frac{P(z)}{Q(z)}$,  where
$P(z)$ and $Q(z)$ are  co-prime polynomials  with degrees   $m$ and $n$,
respectively.  Then  there exists a
constant $c\neq 0$ such that
$$
\lim_{z\rightarrow \infty}R(z)z^{n-m}=c.
$$
As consequence, there exists a constant $M_0>1$ such that
$$
\frac{|c|}{2}|z|^{m-n}\leq |R(z)|\leq 2|c| |z|^{m-n},\ \ \ |z|>M_0.
$$
$ R\in L^p([M_0,\infty)) $ implies that $p(m-n)<-1$, and  so for
$y\geq 0$,
$$
\int_{|x|>M_0}|R(x+iy)|^pdx\leqslant
(2|c|)^p\int_{|x|>M_0}|x+iy|^{p(m-n)}dx
$$
$$
\leqslant (2|c|)^p\int_{|x|>M_0}|x|^{p(m-n)}dx\leqslant
\frac{2^{p+1}|c|^p}{p(n-m)-1}<\infty.
$$
Similarly for $ y>M_0$,
$$
\int_{|x|\leq M_0}|R(x+iy)|^pdx\leqslant (2|c|)^p\int_{|x|\leq
M_0}|x+iy|^{p(m-n)}dx
$$
$$
\leqslant 2(2|c|)^pM_0^{p(m-n)+1}<\infty.
$$
If $R(z)$ is analytic in the upper half  plane $ \mathbb{C}_{+1}$, then
$ Q(z)\neq 0$ for $ z\in \mathbb{C}_{+1}$. If, furthermore,  $ Q(x)\neq
0$ for $ x\in \mathbb{R}$, then $R(z)$ is continuous in the
rectangle  $ E_0=[-M_0,M_0]\times [0,M_0]$, and so $ R \in
H^p(\mathbb{C}_{+1})$.  Otherwise, the null  set  $ N_Q=\{a\in
\mathbb{R}: Q(a)=0\}$ of $Q$ in $\mathbb{R}$ is  a finite set. Let
$N_Q=\{a_1,a_2,\cdots,a_q\}$ with $a_1<a_2<\cdots <a_q$,  and
$P(a_k)\neq 0 (\ k=1,2,\cdots,q)$. Then there exists a
polynomial $Q_1(z)$ with $Q_1(a_k)\neq 0\ (\ k=1,2,\cdots,q)$ and
positive integers $l_k (\  k=1,2,\cdots,q)$ such that
$$
Q(z)=(z-a_1)^{l_1}(z-a_2)^{l_2}\cdots (z-a_q)^{l_q}Q_1(z);
$$
and,
there exist positive constants $\delta , \varepsilon_0 $ and
$M_1>\varepsilon_0$ such that
$$
\varepsilon_0\leq |R(z)(z-a_k)^{l_k}|\leq M_1,
$$
for $ z=x+iy\in I_k=\{z=x+iy:0<|x-a_k|\leq \delta, \ 0\leq y\leq
\delta \}.$

Therefore,
$$
\int_{|x-a_k|\leq \delta}|R(x)|^pdx\geq
\varepsilon_0^p\int_{|x-a_k|\leq \delta}|x-a_k|^{-pl_k}dx.
$$
The fact that $ R\in L^p([a_k-\delta,a_k+\delta]) $ implies that $pl_k<1.$ So,
for $y\in [0,\delta]$,
$$
\int_{|x-a_k|\leq \delta}|R(x+iy)|^pdx\leqslant
M_1^p\int_{|x-a_k|\leq \delta}|x+iy-a_k|^{-pl_k}dx $$
 $$
\leqslant M_1^p\int_{|x-a_k|\leq
\delta}|x-a_k|^{-pl_k}dx=\frac{2M_1^p\delta^{1-pl_k}}{1-pl_k}<\infty.
$$
 Since the poles of $R(z)$ in  the closed upper half plane are identical with  $N_Q$,   $R(z)$ is continuous  in the  bounded closed set
$$\{z\in E_0: z\not\in I_k, \ k=1,2,\cdots,q \}.
$$
 Therefore
$$
\int_{|x|\leq M_0}|R(x+iy)|^pdx
$$
is uniformly bounded for  $y\in [0,M_0]$. This proves that  $ R\in
H^p(\mathbb{C}_{+1})$. If $R(z)$ is analytic in the lower  half  plane $ \mathbb{C}_{-1}$,  Lemma 1 can be  proved similarly.

{\bf  Lemma  2 }\ If  $0<p\leq1$,  $ f\in L^p(\mathbb{R})$ , then, for $\varepsilon>0$,
 there exists a sequence of rational functions $\{R_k(z)\}$, whose  poles are  either $i$ or $-i,$ such that
\begin{eqnarray}\label{cite7}
  \sum_{k=1}^{\infty }||R_k||^p_p\leq (1+\varepsilon )\|f\|_p^p
\end{eqnarray}
and
\begin{eqnarray}\label{cite8}
   \lim_{n\rightarrow \infty}||f-\sum_{k=1}^{n}R_k||_p
  =0.
  \end{eqnarray}
 {\bf Proof }\ \  For the case $0<p<1$, we can assume that $ \|f\|_p^p>0$. The fractional linear mapping (the Cayley Transformation)
$$
z=\alpha (w)=i\frac{1-w}{1+w}
$$
is a conformal mapping from the unit disc $U=\{w:|w|<1\}$ to the
upper half plane  $\mathbb{C}_{+1}$, its inverse mapping is
$$ \beta (z)=\frac{i-z}{z+i}.
$$
Let $x= \alpha (e^{i\theta}),\ \theta \in [-\pi,\pi]$. Then $ x=\tan
\frac{\theta}{2}$ and  $dx=\frac{d \theta}{1+\cos \theta}$. So,
$$
\int_{-\infty}^{\infty}|f(x)|^pdx=\int_{-\pi}^{\pi}\left |f(\tan
\frac{\theta}{2})\right |^p \frac{d \theta}{1+\cos \theta} <\infty.
$$
This implies that  the function
$$g(\theta)=\frac {f(\tan \frac{\theta}{2})
}{(1+\cos \theta)^{\frac{1}{p}}}\in L^p([-\pi,\pi]).
$$
Since the set of trigonometric polynomials is dense in $ L^p([-\pi,\pi])$,
there exists a sequence of rational functions $\{r_k(w)\}$, whose
poles can only be zero,
  with the expression $r_k(e^{i\theta})= \sum_{j=-m_k}^{m_k}c_{k,j}e^{ij\theta}$,
     such that
 $$
   \lim_{k\rightarrow \infty}||g(\theta)-r_k(e^{i\theta})||_{L^p([-\pi,\pi])}
  =0.
$$
Furthermore, for any $ \varepsilon>0$,  the  sequence of rational functions $\{r_k(w)\}$  can
be chosen so that
$$
   ||g(\theta)-r_k(e^{i\theta})||^p_{L^p([-\pi,\pi])}\leq
   \frac{A_\varepsilon }{4^{k+3}},
$$
where $A_\varepsilon =\|f\|_p^p \varepsilon $.
 Since $0<p<1$, there exists a
positive integer $l_p$ such that $ 1<p2^{l_p}\leq 2. $ Take
$m=2^{l_p-1}.$ Then $m$ is a positive integer satisfying $ 1<2pm\leq
2. $ Thus  we have $0\geq 2(pm-1)>-1,$ and, as consequence, the  function
$$
g_1(\theta)=(2\sin^2\theta)^{pm-1}\in L^1\left [0,\frac{\pi}{2}\right ].
$$
The function $ g_2(x)=x^{\frac{1}{p}-m}$ is also continuous in the
interval $ [0,2].$ The Weierstrass Theorem asserts that there exists a sequence of polynomials
$\{q_k(x)\}$ such that
\begin{eqnarray}\label{cite9}
|g_2(x)-q_k(x)|<\frac{A_\varepsilon }{M_k^pC_14^{k+3}},
\end{eqnarray}
where
$$
M_k=\sum_{j=-m_k}^{m_k}|c_{k,j}|+1,\ \ \
C_1=\int_{0}^{\frac{\pi}{2}}g_1(\theta)d\theta.
$$
Thus we obtain
$$
\int_{0}^{\frac{\pi}{2}}|(2\sin^2\theta)^{\frac{1}{p}-m}-q_k(2\sin^2\theta
)|^pg_1(\theta)d\theta\leq \frac{A_\varepsilon}{M_k4^{k+3}}.
$$
The function
$$
s_k(e^{i\theta})=r_k(e^{i\theta})q_k(1+\cos \theta)(1+\cos \theta)^m
$$
is a trigonometric polynomial, and satisfies
$$ J_k=\int_{-\pi}^{\pi}\left|r_k(e^{i\theta})-\frac {s_k(e^{i\theta})}{(1+\cos \theta)^{\frac{1}{p}}}\right |^p d \theta.
$$
$$
\leq M_k^p\int_{-\pi}^{\pi}|1-q_k(1+\cos \theta)(1+\cos
\theta)^{m-\frac{1}{p}}|^p d \theta
$$
$$
 = M_k^p\int_{-\pi}^{\pi}|(1+\cos
\theta)^{\frac{1}{p}-m}-q_k(1+\cos \theta)|^p(1+\cos \theta)^{pm-1}
d \theta$$ $$ =  M_k^p\int_{-\pi}^{\pi}|g_2(1+\cos
\theta)-q_k(1+\cos \theta)|^p(1+\cos \theta)^{pm-1} d \theta.
$$
Hence, by (\ref{cite9}),
$$
 J_k\leq  \frac{A_\varepsilon}{C_14^{k+3}}\int_{-\pi}^{\pi}(1+\cos \theta)^{pm-1} d
\theta= \frac{A_\varepsilon}{C_14^{k+2}}\int_{-\frac{\pi}{2}}^{\frac{\pi}{2}}(2\cos
\theta)^{pm-1} d \theta
$$
$$
\leq
\frac{A_\varepsilon}{C_14^{k+2}}\int_{-\frac{\pi}{2}}^{\frac{\pi}{2}}(2\sin
\theta)^{pm-1} d \theta=\frac{2A_\varepsilon}{4^{k+2}}.
$$
Finally, the function
 $$
  g_k(\theta)=\frac
{s_k(e^{i\theta})}{(1+\cos \theta)^{\frac{1}{p}}}
$$
satisfies
$$ ||g-g_k||^p_{L^p([-\pi,\pi])}
$$
$$
 \leq ||g-r_k(e^{i\cdot})||^p_{L^p([-\pi,\pi])}+||r_k(e^{i\cdot
})-g_k||^p_{L^p([-\pi,\pi])}\leq \frac{A_\varepsilon}{4^{k+1}},
$$
and
$$
||g-g_k||^p_{L^p([-\pi,\pi])}=\int_{-\pi}^{\pi}\left|f(\tan
\frac{\theta}{2})-s_k (e^{i\theta})\right |^p \frac{d \theta}{1+\cos \theta}
$$
$$
= \int_{-\infty}^{\infty}\left |f(x)-s_k\left(\frac{i-x}{x+i}\right
)\right |^pdx \leq \frac{A_\varepsilon}{4^{k+1}}.
$$
The function   $$ Q_k(z)=s_k\left (\frac{i-z}{z+i}\right )
$$
 is a rational
function  whose  poles are either $i$ or  $-i$, and
$$
  ||Q_k||^p_p=\int_{-\infty}^{\infty}|Q_k(x)|^pdx=\int_{-\infty}^{\infty}\left |s_k\left (\frac{i-x}{x+i}\right)\right|^pdx\leq \|f\|_p^p+\frac{A_\varepsilon}{4^{k+1}}
$$
and
$$
   ||f-Q_k||^p_p=\int_{-\infty}^{\infty}\left |f(x)-s_k\left(\frac{i-x}{x+i}\right )\right|^pdx\leq
\frac{2A_\varepsilon}{4^{k+1}} .$$ Therefore,  the  sequence of
rational functions $\{Q_k(z)\}$ can be chosen so that
$$
   ||Q_k-Q_{k-1}||^p_p\leq
   \frac{A_\varepsilon}{4^{k}}.   \ \ \ (k=2,3,\cdots )
$$
Let $$ R_1(z)=Q_1(z), \ R_k(z)=Q_k(z)-Q_{k-1}(z),   \ \ \
(k=2,3,\cdots).
$$
$\{R_k(z)\}$ is  a sequence  of rational functions
whose  poles
 can only be  $i$ or $-i$, satisfying (\ref{cite7}) and (\ref{cite8}).
This completes  the proof of  Lemma 2.

{\bf  Lemma  3 } \ Suppose that   $0<p<1$ and that   $ R\in
L^p(\mathbb{R})$ is  a rational function whose  poles are contained in $\{i,-i\}$, then   there exist two  rational functions   $ P\in H^p(\mathbb{C}_{+1})$ and  $ Q\in
H^p(\mathbb{C}_{-1})$ such that $R(z)=P(z)+Q(z)$ and $$
\|P\|^p_{H^p_{+1}}+\|Q\|^p_{H^p_{-1}}\leq
\left (1+\frac{4\pi}{1-p}\right )\|R\|^p_p,
$$
\\
 {\bf Proof } Let $0<p<1$,   $ R\in L^p(\mathbb{R}),$ and  $R$ be
a rational function whose poles are contained in $\{i,-i\}.$ Then
  $R(z)$ can be written as
 $$
 R(z)=\sum_{k=-n}^{n}c_k( \beta (z))^k,  \ \ \mbox{where }\ \   \beta (z)=\frac{i-z}{z+i}.
$$
Therefore, $\beta (x)=e^{i\theta (x)}$,  where  $\theta (x)=\arg (i-x) -\arg (x+i)\in (-\pi, \pi)$ for $x\in \mathbb{R}$. Define, for each $\varphi \in \mathbb{R}$,
 $$
 P(z,\varphi)=\frac{(\beta (z))^{m}R(z)}{(\beta (z))^m-e^{i\varphi}},\ \ \
 Q(z,\varphi )=\frac{(\beta (z))^{-m}R(z)}{(\beta (z))^{-m}-e^{-i\varphi}},
 $$
  where $m$ is any positive integer greater than the positive integer
  $n$. By Fubini's theorem,
$$
I=\int_{-\pi}^{\pi}\int_{-\infty}^{+\infty}|P(x,\varphi)|^pdxd\varphi=
\int_{-\pi}^{\pi}\int_{-\infty}^{+\infty}\frac{|\beta
(x)|^{mp}|R(x)|^p}{|(\beta (x))^m-e^{i\varphi}|^p}dxd\varphi
$$
$$
=\int_{-\infty}^{+\infty}\int_{-\pi}^{\pi}\frac{|R(x)|^p}{|1-e^{i(\varphi -m\theta (x))}|^p}d\varphi dx.
$$
Observing that
$$\int_{-\pi}^{\pi}\frac{2^pd\varphi}{|1-e^{i\varphi-im\theta
(x)}|^p}=\int_{-\pi}^{\pi}\frac{2^pd\varphi}{|1-e^{i\varphi}|^p}=\int_{-\pi}^{\pi}\frac{d\varphi}{
\sin ^{p}\frac{\varphi}{2}}\leq
4\int_{0}^{\frac{\pi}{2}}\frac{d\varphi}{
(\frac{2}{\pi}\varphi)^p}\leq \frac{2\pi}{1-p},
$$
we obtain that
$$
I\leq \frac{2^{1-p}\pi}{1-p}\int_{-\infty}^{+\infty}|R(x)|^pdx.
$$
Therefore, there is a real number $\varphi$ such that
$$
\int_{-\infty}^{+\infty}|P(x,\varphi)|^pdx\leq
\frac{2\pi}{1-p}\int_{-\infty}^{+\infty}|R(x)|^pdx.
$$
For this specially chosen real number $\varphi$,  by defining
$P(z)=P(z,\varphi),  \ Q(z)= Q(z,\varphi ),
 $
we have $R(z)=P(z)+Q(z)$. Since $m>n $, the  functions $P$ and $Q$  are
rational functions and the poles of
 $P(z)$ and
 $Q(z)$ both are contained in the set $\{i\}\bigcup
 \{x_k:k=0,1,2,\cdots,n-1\}$, where through the Cayley Transformation
 $$
 x_k=\alpha (e^{\frac{i}{n}(\varphi
 +2k\pi)})=\tan^2(\frac{1}{2n}(\varphi
 +2k\pi))$$
are real numbers. Therefore,  $P(z)$ is analytic in the upper half
 plane $\mathbb{C}_{+1}$, and $Q(z)$ is analytic in the lower half plane $\mathbb{C}_{-1},$
 and
$$
\int_{-\infty}^{+\infty}|P(x)|^pdx\leq
\frac{2\pi}{1-p}\int_{-\infty}^{+\infty}|R(x)|^pdx
$$
$$
\int_{-\infty}^{+\infty}|Q(x)|^pdx\leq \left
(1+\frac{2\pi}{1-p}\right )\int_{-\infty}^{+\infty}|R(x)|^pdx.
$$
By Lemma 1, $ P\in H^p(\mathbb{C}_{+1})$, $ Q\in H^p(\mathbb{C}_{-1})$.
This completes  the proof of Lemma 3.\\

{\bf Proof of Theorem 1}\ \ \
   \ If $f\in H^p(\mathbb{C}_{+1}),\ Np>1$, then, for any $\varepsilon>0$,
by Theorem A,  there exists function $f_N$ in  $ H^p(\mathbb{C}_{+1})\bigcap C^\infty (\overline{\mathbb{C}_{+1}})$
such that   $$
\lim\limits_{|z|\rightarrow 0, Im z\geq 0}|z|^{N+1}|f_N(z)=0
$$
 and
 $$\|f_N-f\|_{H^p_{+1}}<\varepsilon.
 $$
The fractional linear mapping (the Cayley Transformation)
$$
z=\alpha (w)=i\frac{1-w}{1+w}
$$
is a conformal mapping from the unit disc $U=\{w:|w|<1\}$ to the
upper half plane  $\mathbb{C}_+$, its inverse mapping is
$$
w=\beta (z)=\frac{i-z}{z+i}.
$$
Let $h_N(w)= f_N(\alpha(w))$ and $ h_N(-1)=0$, then $h_N(w) $ is continuous in the closed disc $\overline{U}$ and
$$
h_N(w)\left(i\frac{1-w}{1+w}\right)^{N+1}\rightarrow 0, \ \ \ w \in \overline{U}\setminus \{-1\} , w\rightarrow -1.
$$
So,
$$
\frac{h_N(w)}{(1+w)^{N+1}} \rightarrow 0, \  \ \ w\rightarrow -1,\ |w|\leq 1,\ w\neq 1.
$$
If let $\tilde{h}_N(w)= \frac{h_N(w)}{(1+w)^{N+1}} $ and $\tilde{h}_N(-1)=0$,
then  $\tilde{h}_N(w)$ is analytic in the unit disc $U$ and continuous in the closed unit disc $\overline{U}$.
Therefore,  there exists polynomial $P_N$ such that
 $$
 \left |\frac{h_N(w)}{(1+w)^{N+1}}-P_N(1+w)\right |< \varepsilon,\ \ |w|\leq 1, w\neq -1.$$
Thus,
$$
|f_N(\alpha(w))-(1+w)^{N+1}P_N(1+w)|< \varepsilon |1+w|^{N+1}, |w|\leq 1, w\neq -1.
$$
Since $z=\alpha(w)$ and $w=\frac{i-z}{i+z}$, the above inequality become
$$
 \left |f_N(z)-\left (\frac{2i}{i+z}\right)^{N+1}P_N\left (\frac{2i}{i+z}\right )\right |< \varepsilon \left |\frac{2i}{i+z}\right |^{N+1},\ Im z\geq 0.$$
Therefore, we obtain
 \begin{equation*}
 \int_{-\infty}^\infty|f_N(x+iy)-R(x+iy)|^p\ dx\leq  \varepsilon^p2^{(N+1)p}\int_{-\infty}^\infty  \left |\frac{1}{x^2+1}\right |^{(N+1)p}dx,
\end{equation*}
where $R(z)= (\frac{2i}{i+z})^{N+1}P_N(\frac{2i}{i+z})\in \mathfrak{R}_N(i)$. This concludes that the class $\mathfrak{R}_N(i)$ is dense in  $ H^p(\mathbb{C}_{+1})$. The proof of Theorem 1 is complete.

The Corollary  can be  proved similarly.

{\bf Proof of Theorem 2}\ \ \
According to Lemma 1 and 2, there exist two  sequences of rational functions $\{P_k(z)\}$ and $\{Q_k(z)\}$ such that  $ P_k\in H^p(\mathbb{C}_{+1})$, $ Q_k\in
H^p(\mathbb{C}_{-1})$,
   $$
\sum_{k=1}^{\infty} \left (\|P_k\|^p_{p}+\|Q_k\|^p_{p}\right )\leq
2\left (1+\frac{2\pi}{1-p}\right )\|f\|^p_p
$$
and
$$
   \lim_{n\rightarrow \infty}||f-\sum_{k=1}^{n}(P_k+Q_k)||_p
  =0.
$$
Since  $$
\|P_k\|^p_{H^p_{+1}}=\|P_k\|^p_{p}\  \ {\rm and} \  \ \|Q_k\|^p_{H^p_{-1}}=\|Q_k\|^p_{p},
$$
we have that (\ref{cite2}) and (\ref{cite3}) hold. For any $\delta>0, y>0$, the functions
$ |P|^p$ and $|Q|^p$ are  subharmonic.  Hence,
$$
  \left |\sum_{k=1}^{n}P_k(x+iy+i\delta)\right |^p\leq  \sum_{k=1}^{n}|P_k(x+iy+i\delta)|^p \leq \frac{2}{\pi \delta } \sum_{k=1}^{n}\|P_k\|^p_{p}.
$$
This implies that  the series
$$
\sum_{k=1}^{\infty}P_k(z)
$$
uniformly converges in the closed upper half plane
$\{z:\mbox{Im}z\geq \delta \}$ for any $\delta>0$. As consequence, the function
$g(z)$ is analytic in the upper half plane $\mathbb{C}_{+1}$.
Similarly, we can prove that the function $h(z)$ is analytic in the
lower half plane $\mathbb{C}_{-1}$.
   (\ref{cite2}) implies that (\ref{cite4}) holds. Therefore,  the non-tangential boundary
limits  $g(x)$ and $ h(x)$  of  functions for $g\in
H^p(\mathbb{C}_{+1})$ and $h\in  H^p(\mathbb{C}_{-1})$ exist almost
everywhere. (\ref{cite3}) implies that $f(x)=g(x)+h(x)$  almost
everywhere.

{\bf A new proof of Theorem B}\ \ \  There exist $f(z)\in H^p(\mathbb{C}_{+1}), \ g(z)\in H^p(\mathbb{C}_{-1})$ such that
$f(x)=g(x), \ \ a.e. x\in \mathbb{R}.$
By Theorem 1 and Corollary 1, for any $\varepsilon >$,   there exist  $R\in \mathfrak{R}_N(i)$ and $R_2\in \mathfrak{R}_N(-i)$
such that
$$
\|f-R_1\|_{H^p_{+1}}=\|f-R_1\|_{p}<\frac{\varepsilon}{4},\ \ \|g-R_2\|_{H^p_{-1}}=\|f-R_2\|_{p}<\frac{\varepsilon}{4}.
$$
By the definition of  $R\in \mathfrak{R}_N(i)$ and $R_2\in \mathfrak{R}_N(-i)$, there exist polynomials $P_1$ and $P_2$
such that
$$
R_1(z)=P_1(\beta(z)+1)(\beta(z)+1)^{N+1}, \ \ R_2(z)=P_2((\beta(z))^{-1}+1)((\beta(z))^{-1}+1)^{N+1},
$$
where $\beta(z)= \frac{i-z}{i+z}. $

  Let $m> \max\{degP_1,\ degP_2\}+N+1$, and define, for each $\varphi \in \mathbb{R}$,
  $$
R(z,  \varphi) = R_1(z) - \frac{(\beta(z))^m(R_1(z)-R_2(z))}{(\beta(z))^m-e^{i\varphi}}.
 $$
Notice that   $\beta (x)=e^{i\theta (x)}$,  where  $\theta (x)=\arg (i-x) -\arg (x+i)\in (-\pi, \pi)$ for $x\in \mathbb{R}$.
  By Fubini's theorem,
  \begin{equation*}
 \begin{split}
& J=\int_{-\pi}^{+\pi}\int_{-\infty}^{+\infty}|R(x, \varphi)-R_1(x)|^p\ dx d\varphi \\
&= \int_{-\infty}^{+\infty}\int_{-\pi}^{+\pi}\frac{|R_1(x)-R_2(x)|^p}{|(1-e^{i\varphi-im \theta(x)}|^p}\  d\varphi dx.\\
\end{split}
\end{equation*}
 Observing
$$\int_{-\pi}^{\pi}\frac{2^pd\varphi}{|1-e^{i\varphi-im\theta
(x)}|^p}=\int_{-\pi}^{\pi}\frac{2^pd\varphi}{|1-e^{i\varphi}|^p}=\int_{-\pi}^{\pi}\frac{d\varphi}{
\sin ^{p}\frac{\varphi}{2}}\leq
4\int_{0}^{\frac{\pi}{2}}\frac{d\varphi}{
(\frac{2}{\pi}\varphi)^p}\leq \frac{2\pi}{1-p},
$$
we obtain
$$
J\leq \frac{2^{1-p}\pi}{1-p}\int_{-\infty}^{+\infty}|R_1(x)-R_2(x)|^pdx.
$$
Therefore, there is a real number $\varphi$ such that
 \begin{equation*}
\int_{-\infty}^{+\infty}|R(x, \varphi)-R_1(x)|^p\ dx \leq  \frac{2\pi}{1-p}((\varepsilon/4)^p+(\varepsilon/4)^p).
\end{equation*}
Therefore, we have
 \begin{equation*}
 \begin{split}
&\int_{-\infty}^{+\infty}|R(x, \varphi)-f(x)|^p\ dx  \\
& \leq \int_{-\infty}^{+\infty}|R(x, \varphi)-R_1(x)|^p\ dx +\int_{-\infty}^{+\infty}|R_1(x)-f(x)|^p\ dx\\
& \leq  (\varepsilon/4)^p+\frac{4\pi}{1-p}(\varepsilon/4)^p.
\end{split}
\end{equation*}
So, $R(z)=R(z,\varphi)\in L^p(\mathbb{R})$ is a rational function of $z.$ There is a polynomial $P_3$ with ${\rm deg} P_3=N+1+{\rm deg} P_1$
such that $ R(z)=P_3( \beta (z)+1)$. So the poles of $R$ are contained in $\{x_k:\ k= 0,1,...,m+1\}$, where
$$x_k= \alpha (e^{\frac{i(\varphi+2k \pi)}{m}})=\tan^2(\frac{(\varphi+2k \pi)}{2^m}).$$
Thus, $R(z)\in X^p$.

{\bf Proof of   Theorem 3}.
 Recall that
  the Paley-Wiener Theorem asserts that
  $g\in
H^2(\mathbb{C}_{+1})$  if and only if  $\hat{g}\in L^2(\mathbb{R})$ with the support supp$\hat{g}\subset[0,\infty)$,
such that
$$g(z)=\frac{1}{\sqrt{2\pi}}\int^\infty_0 \hat{g}(t)e^{itz} dt \quad\quad (z\in \mathbb{C}_{+1}).
$$
In the case there holds the equality
  $$ \int^\infty_0 |\hat{g}(t)|^2dt=\|g\|_{H_{+1}^2}^2.
  $$
 Let  $0<p\le 1,\ f\in H^p(\mathbb{C}_{+1}).$ For $\delta>0$, let
$f_\delta(z)=f(z+i\delta)$. Then  $ |f|^p$ is subharmonic, and, for
$y>0$,
$$
|f_\delta(x+iy)|\leqslant C_p\|f\|_{H_{+1}^p}\delta^{-\frac{1}{p}},
$$
 where  $C^p_p=\frac{2}{\pi}$. Therefore
$$
\int_{-\infty}^{\infty}|f_\delta(x+iy)|^2dx\leqslant
\int_{-\infty}^{\infty}|f_\delta(x+iy)|^p|f_\delta(x+iy)|^{2-p}dx\leqslant
C_p^{2-p}\|f\|_{H_{+1}^p}^2 \delta^{1-\frac{2}{p}},
$$
and
$$
\int_{-\infty}^{\infty}|f_\delta(x+iy)|dx=
\int_{-\infty}^{\infty}|f_\delta(x+iy)|^p|f_\delta(x+iy)|^{1-p}dx\leqslant
C_p^{1-p}\|f\|_{H_{+1}^p} \delta^{1-\frac{1}{p}}.
$$
 Therefore,  ${\rm supp}\hat{f_\delta}\subset[0,\infty),$ and
\begin{eqnarray}\label{cite10}
f_\delta(z)=\frac{1}{\sqrt{2\pi}} \int_{0}^{\infty}{\hat
f}_\delta(s)e^{itz} dt.
\end{eqnarray}
For $y>0$, $f_\delta(x+iy)=(P_y\ast f_\delta)(x)$, where
$$
P_y(x)=\mbox {Re}\left(\frac{i}{\pi z}\right)= \frac{y}{\pi
(x^2+y^2)}
$$
is the Poisson  kernel of the upper plane $\mathbb{C}_+$.
  It is well known that  $ f_\delta \in L^2(\mathbb{R}),\ P_y\in L^1(\mathbb{R}),$
$\hat P_y(s)=e^{-|s|y}$ for almost all $s\in \mathbb{R}$, and  $ \hat
f_{\delta+y}(s)= \hat f_{\delta}(s)e^{-|s|y}$. So, for almost all $s\in
\mathbb{R}$,  $ \hat f_{\delta+y}(s)e^{|s|(\delta+y)}= \hat
f_{\delta}(s)e^{|s|\delta}$. Hence, the function  $F(s)= \hat
f_{\delta}(s)e^{|s|\delta}$ is independent of $\delta>0$, with
 supp$F\subset[0,\infty),$ and
$$
\int^\infty_{-\infty}|F(s)|^2
e^{-2|s|\delta}dt=\int^\infty_{-\infty}|\hat f_\delta(x)|^2 dx
$$
$$
= \int^\infty_{-\infty}|f_\delta(x)|^2 dx
  \leqslant  C_p^{2-p}\|f\|_{H_{+1}^p}^2 \delta^{1-\frac{2}{p}}.
$$
Therefore, for any $\delta>0$,
 $$
|F(s)|=| \hat f_{\delta}(s)|e^{s\delta}\leqslant
\|f_{\delta}\|_1e^{s\delta} \leqslant
C_p^{1-p}\|f\|_{H_{+1}^p}e^{s\delta} \delta^{-B_p},
$$
where $B_p=\frac{1}{p}-1\geqslant 0$. Since
$$\inf\{ |s|\delta-B_p\log
\delta:\delta>0\}=B_p-B_p(\log B_p-\log |s|),
$$
 we have
$$
|F(s)| \leqslant C_p^{1-p}\|f\|_{H_{+1}^p}B_p^{-B_p}e^{B_p}|s|^{B_p}.
$$
 Thus   $F$ is a slowly increasing  continuous function $F$ whose  support is contained in $[0,\infty)$.                                                                                                                                                                                                                                                                                                                                                                                                                                                                                                                                                                                                                                                                                                                                                                                                                                                                                                                                                                                                                                                                                                                                                                                                                                                                                                                                                                                                                                                                                                                                                                                                                                                                                                                                                                                                                                                                                                                                                                                                                                                                                                                                                                                                                                                                                                                                                                                                                                                                                                                                                                                                                                                                                                                                                                                                                                                                                                                                                                                                                                                                                                                                                                                                                                                                                                                                                                                                                                                                                                                                                                                                                                                                                                                                                                                                                                                                                                                                                                                                                                                                                                                                                                                                                                                                                                                                                                                                                                                                                Letting $\delta \rightarrow 0$ in (\ref{cite10}), we see that (\ref{cite7}) holds.  $F$ can also be regard as a tempered distribution defined through
 $$
(F,\hat{\varphi})=\int_{\mathbb{R}}F(x)\hat{\varphi}(x)dx \ \
$$
 for $\varphi$ in the Schwarz class $\mathbb{S}$. So,
 $$
\lim_{\delta\rightarrow 0}\int_{\mathbb{R}}f_\delta(x)\varphi(x)dx
=\lim_{\delta\rightarrow0}\int_{0}^{+\infty}\hat{f}_\delta
(x)\hat{\varphi}(x)dx
$$
$$
=\lim_{\delta\rightarrow0}\int_{0}^{+\infty}F(x)e^{-\delta x}\hat{\varphi}(x)dx=(F,\hat{\varphi}). \ \
$$
 This completes the proof of Theorem 3.

{\bf A proof of Theorem C}.
 Let $1\leq p\leq 2$. If  $f\in L^{p}$  and ${\rm supp} \hat{f}\subset [0,\infty)$, then
$$
|\chi_{[0,\infty)}(t)e^{2\pi iz\cdot t}\hat{f}(t)| = \chi_{[0,\infty)}(t)|\hat{f}(t)|e^{-2\pi y\cdot t} \in L^{1}(\mathbb{R}^{n}),
$$
 where $\chi_{[0,\infty)}(t)$ is the characteristic function of $[0,\infty)$, that is, $\chi_{[0,\infty)}(t)= 1,$
for $t \in [0,\infty) $, and otherwise zero.
It is evident that the function
$$
G(z)=\frac{1}{\sqrt{2 \pi}}\int_{\mathbb{R}}e^{ izt}\hat{f}(t) dt=\int_{\mathbb{R}}
 \chi_{[0,\infty)}(t)e^{ iz t}\hat{f}(t) dt
$$
is holomorphic in $\mathbb{C}_{+1}$.
To complete the proof of Theorem C, it is sufficient to prove that $G(z)\in H^{p}(\mathbb{C}_{+1})$ and the boundary limit of $G(z)$ is $ f(x)$ as $ y \rightarrow 0$. Fix $z\in \mathbb{C}_{+1}$ and let
 $$
 g_z(t)=\chi_{[0,\infty)}(t)\frac{e^{izt}}{\sqrt{2 \pi}},\ \  \widetilde{g}_z(t)=g_z(-t)
 $$
 then $g_z\in L^1(\mathbb{R})\bigcap L^2(\mathbb{R})$,
  $
 \hat{g}_z(s)=\frac{1}{2\pi i(s-z)}
 $ and
 $$
\begin{array}{rl}
G(z)=& \frac{1}{\sqrt{2 \pi}}  \int_{\mathbb{R}} \chi_{[0,\infty)}(t)e^{ iz t}\frac{1}{\sqrt{2 \pi}}(\int_{\mathbb{R}}e^{- is t}F(s)ds) dt\\
=& \frac{1}{2\pi i}\int_{\mathbb{R}} \frac{f(s)ds}{s-z}.
\end{array}           .
$$
 For $z,w\in \mathbb{C}_{+1}$, let
$$
 I(z,w)=\frac{1}{4 \pi ^2}\int_{\mathbb{R}}\frac{f(t)dt}{(t-z)(t-w)}.
$$
 Then
 $$
 I(z,w)=\int_{\mathbb{R}}\hat{g}_z(t)f(t)\hat{g}_w(-t) dt.
$$
For $z,w\in \mathbb{C}_{+1}$, $\sqrt{2\pi}\hat{g}_z(t)\hat{g}_w(-t)=\widehat{g_z*\widetilde{g}_w}(t)$,
 where   $$
(g_z*\widetilde{g}_w)(t)=\int_{\mathbb{R}}g_z(\xi ) \widetilde{g}_w(t-\xi)d\xi =\int_{\mathbb{R}}g_z(\xi ) g_w(\xi-t)d\xi
$$
$$
=\frac{1}{2\pi}\int_{\mathbb{R}} \chi_{[0,\infty)}(\xi )e^{2\pi iz   \xi } \chi_{[0,\infty)}(\xi-t)e^{2\pi i w (\xi- t)} d\xi.
$$
 Therefore,
$$
I(z,w)=\frac{1}{\sqrt{2\pi}}\int_{\mathbb{R}}\hat{f}(s)\chi_{[0,\infty)}(s)(g_z*\widetilde{g}_w)(s) ds
$$
$$
=\frac{1}{(\sqrt{2\pi})^3}\int_{\mathbb{R}}\hat{f}(s)\chi_{[0,\infty)}(s)\int_{\mathbb{R}}\chi_{[0,\infty)}(\xi )e^{2\pi iz  \cdot \xi } \chi_{[0,\infty)} (\xi-s)e^{2\pi i w \cdot (\xi- s)} d\xi ds.
$$
By Fubini's theorem and the relation
$$
\chi_{[0,\infty)}(t)\chi_{[0,\infty)}(t+s)\chi_{[0,\infty)}(s)=\chi_{[0,\infty)}(t)
\chi_{[0,\infty)}(s),
$$
we have
 $$
\begin{array}{rl}
I(z,w)=& \frac{1}{(\sqrt{2\pi})^3}\int_{\mathbb{R}}\int_{\mathbb{R}}\chi_{[0,\infty)}(s)\chi_{[0,\infty)}(t)\chi_{[0,\infty)}(t+s)e^{iz (s+t)}e^{ i w t}\hat{f}(s )ds dt\\
=&\frac{1}{(\sqrt{2\pi})^3}\int_{\mathbb{R}}\int_{\mathbb{R}} \chi_{[0,\infty)}(t)\chi_{[0,\infty)}(s)e^{iz  s}e^{i(z+ w) t}\hat{f}(s )dt  ds\\
=& \frac{i}{2\pi }\frac{G(z)}{z+w}.
\end{array}
$$
Thus,
 for $z\in \mathbb{C}_{+1}$, we have $-\bar{z}\in \mathbb{C}_{+1}$, and
$$
I(z,-\bar{z})=\frac{i}{2\pi }\frac{G(z)}{z-\bar{z}}=\frac{G(z)}{4\pi y},\ z=x+iy, y>0.
$$
So,
$$
G(z)=\int_{\mathbb{R}}\frac{4\pi y f(t)dt}{(2\pi
)^2(t-z)(t-\bar{z})}=\int_{\mathbb{R}}f(t)P(x-t,y) dt,
$$
where $P(x,y)=\frac{y}{\pi (x^2+y^2)}$ is the Poisson Kernel of the upper half plane  $\mathbb{C}_{+1}$.
Therefore, the boundary limit of $G(z)$ is $ f(x)$ as $ y \rightarrow 0$ and $ G(z) \in H^{p}(\mathbb{C}_{+1}).$
The proof of Theorem C is complete.


\begin{thebibliography}{BSYM}

\bibitem{Ale} A.B. Aleksandrov, { \it Approximation by rational functions and an analogy of the M. Riesz theorem
on conjugate fucntions for $L^p$ with $p\in (0,1)$,} Math.USSR Sbornik, {\bf 35} (1979),301 - 316.

\bibitem{Cim} J.A, Cima and W.T. Ross, {\it The Backward Shift on the Hardy Space}   Mathematical Surveys and  Monographs Vol 79,
American Mathematical Society, Providence, Rhode Island, 2000.

\bibitem{Den} G.T. Deng, {\it Complex Analysis} (in Chinese),
Beijing Normal University Press, 2010.

\bibitem{Du} P. Duren, {\it Theory of $H^p$ Spaces,} Dover Publications. Inc. 2000.

\bibitem{Ga} J.B. Garnett, {\it Bounded Analytic Functions,} Academic Press,New York, 1981.
American Mathematical Society, Providence, Rhode Island, 1996.

\bibitem{Koo} P. Koosis, {\it Introduction to  $H_p$ Spaces,} 2nd ed., Cambridge University Press, 1998.

\bibitem{Le} B. Ya. Levin, {\it Lectures on Entire Functions,} Translations of Mathematical Monographs Vol 150,
American Mathematical Society, Providence, Rhode Island, 1996.

\bibitem{Q1} T. Qian, {\it Characterization of boundary values of functions in Hardy spaces
with application in signal analysis,} Journal of Integral Equations and Applications, Volume {\bf 17} Issue 2 (2005) 159-198.

\bibitem{Q2} T. Qian, Y. S. Xu, D. Y. Yan, L. X. Yan and B. Yu, {\it Fourier Spectrum Characterization of Hardy Spaces and Applications,} Proceedings of the American Mathematical Society, Volume 137, Number 3, March 2009, page 971-980. DOI:10.1090/S0002-9939-08-09544-0.

\bibitem{QW} T. Qian and Y. B. Wang, {\it Adaptive Decomposition Into Basic Signals of Non-negative Instantaneous Frequencies - A Variation and Realization of Greedy Algorithm, } Adv. Comput. Math. 34 (2011), no. 3, 279¡V293.

\bibitem{Ru} W. Rudin, {\it Real and Complex Analysis,} 3rd Edn, New York: McGRAW-HILL International Editions, 1987.

 \bibitem{Wa} J. L. Walsh, {\it Interpolation and Approximation by Rational Functions
in the Complex Plane,} AMS, 1969.

\end{thebibliography}
\end{document}